\documentclass[12pt]{article}
\usepackage{graphicx, amsmath,amsfonts,amsthm, euscript, cite}
\usepackage{verbatim}

\newtheorem{theorem}{Theorem}[section]

\newtheorem{lemma}{Lemma}[section]
\newtheorem{corollary}{Corollary}[section]
\newtheorem{remark}{Remark}[section]
\newtheorem{example}{Example}[section]

\usepackage{graphicx, amsmath}
\usepackage{color}
\definecolor{Blue}{rgb}{0.3,0.3,0.9}

\newtheorem*{Ack}{Acknowledgment}

\begin{document}
\normalsize
\title{Open quantum random walks  with decoherence on coins with $n$ degrees of freedom}

\maketitle
Sheng Xiong
\vskip 3pt
Department of Mathematics\\
        University of Pittsburgh,
        Pittsburgh, PA 15260
\vskip 6pt
  Wei-Shih Yang
\vskip 3pt
Department of Mathematics\\
        Temple University,
        Philadelphia, PA 19122

\vskip 6pt
Email: sxiong@pitt.edu,  yang@temple.edu



\vskip12pt

\begin{abstract}
In this paper, we define a new type of decoherent quantum random walks with parameter $0\le p\le 1$, which becomes a unitary quantum random walk (UQRW) when $p=0$ and an open quantum random walk (OPRW) when $p=1$ respectively.  We call this process  a partially open quantum random walk (POQRW). We study the limiting distribution of a POQRW on $Z^1$ subject to decoherence on coins with $n$ degrees of freedom, which  converges to a convex combination of normal distributions if the superoperator $\mathcal{L}_{kk}$ satisfies the eigenvalue condition, that is, 1 is an eigenvalue of $\mathcal{L}_{kk}$  with multiplicity one and all other eigenvalues have absolute values less than 1.  A   Perron-Frobenius type of theorem is provided in determining whether or not  the superoperator satisfies the eigenvalue condition. Moreover, we compute the limiting distributions of characteristic equations of the position probability functions for $n=2$ and 3.

\end{abstract}

\section{Introduction}

It has been given rise to a vast  field of exploration for the behavior of quantum open systems \cite{APS1, AGS, PS, PS1,KY} since  S. Attal et al. \cite{APS} recently introduced a new type of an open quantum random walk (OQRW) on graphs, which is an exact quantum analogue of classical  Markov  Chain. This new type of walks has exhibited many interesting phenomenon since it is partially a quantum walk and partially a classical walk; moreover,  it has a speed-up property like a quantum walk behavior and yet it has a Markov property like a classical walk \cite{AGS}. In particular,  S. Attal et al.  \cite{APS} pointed out that there is a strong link between OQRWs and the well known unitary quantum random walks (UQRW) and that it is difficult to produce  the limit distribution for OQRWs,  due to lack of knowledge about the invariant measures of this Markov chain and even their existence at time of their paper was written. But soon in  \cite {AGS} S. Attal et al.  obtained  a Central Limit Theorem for the case of the nearest neighbors homogeneous  OQRW on $Z^d$ under an assumption that the superoperator  admits a unique invariant state. The variance is somewhat abstract, so Konno and Yoo \cite{KY} further studied the limit distributions of OQRWs on one-dimensional lattice space and compute the distribution of the OQRWs concretely for many examples and thereby obtain the limit distributions of them.  Under an equivalent condition that  1 is an eigenvalue of the superoperator  with multiplicity one and all other eigenvalues have absolute values less than 1(for convenience, throughout the paper, this condition is called``eigenvalue condition"),   the authors of this paper also  \cite{FFXY}  obtained  the limit distribution (a convex combination of Gaussian distribution) for UQRWs on $Z$. This suggests us to further explore the strong link between these two types of quantum walks and the limit distributions of OQRWs under the eigenvalue condition.

In this paper, we consider  an OQRW  resulting from a total decoherent on the coin space of a UQRW while  the dimension of the coin space can be arbitrarily large.  In such case, we investigate the transitions of how a unitary QRW is eventually collapsed into an OQRW.  Moreover, we study a type of QRWs where each step has a probability $p$ of decoherent on the coin space. If $p>0$, this process will eventually collapse into OQRW and exhibit the diffusive behavior. We call this process partially open quantum walk with parameter $p$ (POQRW). Here $0 \le p \le 1$. If $p=0$, it's a UQRW; if $p=1$, it is corresponding to an OQRW.  Assuming the superoperator of the POQRW satisfies the eigenvalue condition, we conclude that  the POQRW also  converges to a convex combination of normal distributions.

The rest of our paper is organized as follow: in Section 2 we briefly summarize the formalism of POQRWs. In Section 3 we prove a  limiting theorem of  POQRWs on $Z^1$ with decoherence on coins with $n$ degrees of freedom with the help of the Quantum Fourier Transform and the generalized Gell-Mann matrices basis. In Section 4  we prove a Perron-Frobenius type of theorem which is very useful in determining when the superoperator $\mathcal{L}_{kk}$ satisfies the eigenvalue condition. In Section 5 we demonstrate some examples of POQRWs ($n=2, 3$) that satisfy the eigenvalue condition and give  explicit formulas for their limiting distributions.  In particular, our results generalize the Central Limit Theorem in \cite{AGS}.

\setcounter{equation}{0}

\section{The partially open quantum random walk on the lattices $Z^d$ and decoherence}

Let us consider a general open quantum random walk on $Z^d$.   Let $\{e_j\}_{j=1}^{d}$ be the standard orthonormal basis of $Z^d$ and we put  $e_{j+d}= -e_j,  j=1,2,\dots, d$. 
 We denote the state space by a Hilbert space $\mathcal{H}=\mathcal{H}_{P}\otimes \mathcal{H}_{C}$, where  $\mathcal{H}_{P}$ denotes the position space and  $\mathcal{H}_{C}$ denotes the coin space. 
 The orthonormal basis of the position space $\mathcal{H}_{P}$ are $|x>$, where $x\in Z^d$, the basis of the coin space $\mathcal{H}_{C}$ are $|\xi_i>, i=1,2,\dots,n.$  We will assume that the walk starts at the origin.  Let us describe the dynamics of the quantum walker.

Let $U=[u_1,u_2,\dots,u_n]^T$ be a $n \times n$ unitary matrix and $\Pi_j, j=1,2,\dots, m$ be orthogonal projection matrices which partition the matrix $U$ into $m$ matrices $B_j, j=1,2,\cdots, m$,
 where 
\[B_1=[u_1,u_2,\dots,u_{d_1}, 0,0,\dots,0]^T,\] 
 \[ B_2=[0,0,\dots, 0,u_{1+d_1},u_{2+d_1},\dots,u_{d_1+d_2}, 0,0,\dots,0]^T,\] 
 \[B_j=[0,0,\dots, 0,u_{1+\sum_{i=1}^{j-1} d_i},u_{2+\sum_{i=1}^{j-1} d_i},\dots,u_{\sum_{i=1}^{j} d_i}, 0,0,\dots,0]^T,\] 
  for $j=1,2\dots,m,$ and $\sum_{i=1}^{m} d_i=n$.
 Note that 
 \[B_j=\Pi_j U,\; \sum_{j=1}^{m}  B^*_j B_j=I,  \; \sum_{j=1}^{m} B_j=U.\]

Let us define $$L_{x}^{j}=  |x+e_j> <x| \otimes B_j,$$   
$$L^{j}= \sum_{x} L_{x}^{j}, \hskip 12pt \mbox { for } j=1,2,\dots, m,$$  where $m=2d$, and we also put  
 $L_x= \sum_{j=1}^{m} L_{x}^{j},$ and $L=\sum_{x} L_{x}.$ 
 Then $L$ is a unitary operator on $\mathcal{H}$. Let $\psi_0 \in \mathcal{H}$ and $\psi_n=L^n \psi_0.$ 
Then $\{\psi_n\}_{n=0}^{\infty}$ is called a quantum random walk on $Z^d$. 


 In terms of density operator $$M_0(\rho)=L \rho L^*$$ is the one-step dynamics of the density operator.  
So $\rho ^{(n)}=M_0^n (\rho ^{(0)})$ is the dynamics of the quantum random walks.  
Note that 
$$M_0(\rho)= \sum_{j, j^\prime =1}^{m}\sum_{x, x^\prime} { L_{x}^j} \rho { L_{x^\prime}^{j^\prime}}^*.$$
Let $$M_{open}(\rho)= \sum_{j=1}^{m}\sum_{x} { L_{x}^j} \rho { L_{x}^j}^*.$$ 
 Then $\rho ^{(n)}=M_{open}^n (\rho ^{(0)})$ is the dynamics of the open quantum random walks defined by S. Attal, F. Petruccione and I. Sinayskiy \cite{APS} (2012).  
 
  Let us now define the partially open quantum random walk. 
 Let
  \[ \quad D_j=\sqrt {p} L^{j},  j=1,2,\dots, m,\] and \[  \quad D_0=\sqrt{q} L,\] 
 where $0\le p, q\le 1, p+q=1.$  
 We have  $ \sum_{j=0}^{m}{ D_j}^* D_{j}=I_\mathcal{H}$ and therefore  $\{D_j\}_{j=0}^{m}$ can be viewed as measurements  on $\mathcal{H}=\mathcal{H}_{P}\otimes \mathcal{H}_{C}$.
Let us define a positive map  on $ \mathcal{H}$:
  $$M(\rho)= M_p(\rho)=\sum_{j=0}^{m}{ D_{j}} \rho D_{j}^*,$$
 where $\rho$ is a density operator on $\mathcal{H},$ then 
  $\{M^n\}_{n=1}^{\infty}$ is called the POQRW on $Z^d$ associated with $\{B_j\}$. Note that if $p=0$, it is a UQRW. If $p=1$, it is corresponding to an OQRW, as seen below.


The superoperator of open quantum random walks defined in S. Attal, N. Guillotin-Plantard and C. Sabot \cite {AGS} is 
$$M_{open}(\rho)= \sum_{j=1}^{m}\sum_{x} { L_{x}^j} \rho { L_{x}^j}^*,$$ 
and in  our case, $$M_{1}(\rho)= \sum_{j=1}^{m}\sum_{x, x^\prime} { L_{x}^j} \rho { L_{x^\prime}^j}^*,$$
so they are not exactly the same operators. 
But if $\rho =\sum_{x} |x><x |\otimes \rho_x$, where $\rho_x$ is a positive operator on  $\mathcal{H}_{C}$ such that $\sum_x Tr (\rho _x)=1$, then $\rho^1=M_{open}(\rho)$ also has the same form and 
$$M_{open}(\rho)=M_1(\rho), \hskip 12pt \mbox{ for all } \rho =\sum_{x} |x><x |\otimes \rho_x.$$ 
Therefore, our case $p=1$ includes the corresponding open quantum random walks.  
Since we are dealing with partially open quantum random walks, the form of density operator has a nontrivial off-diagonal (in $x$-space) entry after the initial step, we need to consider general density operators for later iterations.  
Thus our formula $M_1$ is a generalization of the open quantum random walks.

\

In general, the density operator for quantum random walk in Fourier transformation basis is given by
\begin{equation}
\rho=\int\frac{dk}{(2\pi)^d}\int\frac{dk^{\prime}}{(2\pi)^d}|k><k^{\prime}|\otimes \chi_{kk^{\prime}},
\end{equation}
where $k=(k_1, k_2, ..., k_d)$, $0 \le k_i \le 2 \pi$, and $dk=dk_1dk_2...dk_d$. 
Then after one step the density operator  becomes
\begin{eqnarray*}
\rho\rightarrow\rho^{\prime}
&=&\sum_{j=0}^{m} D_j \rho D^{*}_j \\
&=&\sum_{j=1}^{m}\sum_{x^{\prime}}\sum_x p( {L_{x}^{j}} \rho {L_{x'}^{j}}^*)+q\sum_{x^{\prime}}\sum_x ( {L_{x}}\rho {L_{x'}}^*) \\
&=&\int\frac{dk}{(2\pi)^d}\int\frac{dk^{\prime}}{(2\pi)^d}|k><k^{\prime}|\otimes \sum_{j=0}^{m}A_jU_k\chi_{kk^{\prime}}U^*_{k^{\prime}}A^*_j,
\end{eqnarray*}
where  $A_j=\sqrt{p} \Pi_j, A_0=\sqrt{q} I, B_{jk}=e^{-ike_j}B_j, U_{k}=\sum_{j=1}^{m} B_{jk}$, and $ke_j$ denotes the dot product of $k$ and $e_j$.  
Suppose the quantum walk starts at a pure state $|0>\otimes |\Phi_0>$,
then the initial state is
\begin{equation}
\rho_0=\int\frac{dk}{(2\pi)^d}\int\frac{dk^{\prime}}{(2\pi)^d}|k><k^{\prime}|\otimes |\Phi_0><\Phi_0|.
\end{equation} 
Let the quantum random walk proceed for $t$ steps. Then the state evolves to
\begin{eqnarray*}
\rho_t=\int\frac{dk}{(2\pi)^d}\int\frac{dk^{\prime}}{(2\pi)^d}|k><k^{\prime}|\otimes \\
\sum_{j_1,\ldots,j_t}A_{j_t}U_k\cdots A_{j_1}U_k |\Phi_0><\Phi_0| U^*_{k^{\prime}}A^*_{j_1}\cdots U^*_{k^{\prime}}A^*_{j_t}.
\end{eqnarray*}

Let  $\mathcal{L}_{kk^{\prime}}$ be the operator acting on the vector space of linear operators $L(\mathcal{H}_C)$ and $\chi_{kk^{\prime}}\in L(\mathcal{H}_C)$, then it follows from the direct computation 
\begin{equation}
\mathcal{L}_{kk^{\prime}}\chi_{kk^{\prime}}\equiv\sum_{j}A_{j}U_k \chi_{kk^{\prime}} U^*_{k^{\prime}}A^*_{j} =p\sum_{j=1}^{m} B_{jk}\chi_{kk^{\prime}}B^*_{jk^{\prime}}+qU_k\chi_{kk^{\prime}}{U_{k^{\prime}}}^*.
\end{equation}
$\mathcal{L}_{kk^{\prime}}$ is a linear operator  that maps from  $L(\mathcal{H}_C)$ to $L(\mathcal{H}_C)$, that is $\mathcal{L}_{kk^{\prime}}\in L(L(\mathcal{H}_C))$. $\mathcal{L}_{kk^{\prime}}$  is also called a superoperator, and  
\begin{eqnarray*}
\int\frac{dk}{(2\pi)^d}\int\frac{dk^{\prime}}{(2\pi)^d}|k><k^{\prime}|\otimes\sum_{j} A_{j}U_k\chi_{kk^{\prime}}  U^*_{k^{\prime}}A^*_{j} \\
=\int\frac{dk}{(2\pi)^d}\int\frac{dk^{\prime}}{(2\pi)^d}|k><k^{\prime}|\otimes \mathcal{L}_{kk^{\prime}}\chi_{kk^{\prime}}.
\end{eqnarray*}
In terms of the superoperator $\mathcal{L}_{kk^{\prime}}$,
\begin{equation}
\rho_t=\int\frac{dk}{(2\pi)^d}\int\frac{dk^{\prime}}{(2\pi)^d}|k><k^{\prime}|\otimes \mathcal{L}^t_{kk^{\prime}}|\Phi_0><\Phi_0|.
\end{equation} 
The probability to reach a point $x$ at time $t$ is
\begin{eqnarray}
p(x,t)&=&Tr\{[ |x><x|\otimes I ]\rho_t\}\nonumber\\
&=&\frac{1}{(2\pi)^{2d}}\int dk\int dk^{\prime}<k|x><x|k^{\prime}>Tr\{ \mathcal{L}^t_{kk^{\prime}}|\Phi_0><\Phi_0|\}\nonumber\\
&=&\frac{1}{(2\pi)^{d}}\int dk\int dk^{\prime} e^{-ix(k-k^{\prime})}Tr\{ \mathcal{L}^t_{kk^{\prime}}|\Phi_0><\Phi_0|\}.
\end{eqnarray}



\section{The limiting distributions of  quantum walks on $Z^1$ with decoherence on coins with $n$ degrees of freedom}
\setcounter{equation}{0}
Consider the Fourier transformation of $p(x,t)$
\begin{equation}
\hat{P}(\nu,t)\equiv<e^{i\nu x}>_t=\sum_xe^{i\nu x}p(x,t).
\end{equation}
To simplify $<e^{i\nu x}>_t$, we use the properties of   $\delta$ function
\begin{equation}
\frac{1}{2\pi}\sum_x x^m e^{-ix(k-k^{\prime})}=(-i)^m\delta^{(m)}(k-k^{\prime}).
\end{equation}
Then
\begin{eqnarray}
<e^{i\nu x}>_t &=& \sum_x e^{i\nu x}p(x,t)\nonumber\\
&=&\sum_x e^{i\nu x}\int \frac{dk}{2\pi}\int \frac{dk^{\prime}}{2\pi} e^{-ix(k-k^{\prime})}Tr\{ \mathcal{L}^t_{kk^{\prime}}|\Phi_0><\Phi_0|\}\nonumber\\
&=&\int \frac{dk}{2\pi}\int \frac{dk^{\prime}}{2\pi}\sum_x e^{ix(\nu+k^{\prime}-k)}Tr\{ \mathcal{L}^t_{kk^{\prime}}|\Phi_0><\Phi_0|\}\nonumber\\
&=&\int \frac{dk}{2\pi}\int \frac{dk^{\prime}}{2\pi}2\pi \delta(\nu+k^{\prime}-k)Tr\{ \mathcal{L}^t_{kk^{\prime}}|\Phi_0><\Phi_0|\}\nonumber\\
&=&\frac{1}{2\pi}\int dk Tr\{ \mathcal{L}^t_{k,k+\nu}|\Phi_0><\Phi_0|\}.
\end{eqnarray}
Let $\hat{O}$ denote any operator on $\mathcal{H}_C$.  Then the generating function of $<e^{i\nu x}>_t$ is given by
\begin{eqnarray*}
G(z,\nu)&=&\sum_{t=0}^{\infty}z^t<e^{i\nu x}>_t\\
&=&\frac{1}{2\pi}\int dk\sum_{t=0}^{\infty}z^tTr\{\mathcal{L}^t_{k,k+\nu}\hat{O}\}\\
&=&\frac{1}{2\pi}\int dk Tr\{\frac{1}{I-z\mathcal{L}_{k,k+\nu}}\hat{O}\}.
\end{eqnarray*}
where $|z|<1$.
 This definition makes sense since the spectrum of $\mathcal{L}_{k,k+\nu}$ is less than or equal to $1$ by Lemma (3.1) in \cite{FFXY}.

Let the coin Hilbert space $\mathcal{H}_C$ be $n$-dimensional spanned by orthogonal basis $\xi_1, \xi_2,\dots, \xi_n $,  $\Pi_1$ be the orthogonal projection onto the subspace spanned by $\xi_1, \xi_2,\dots, \xi_{n_1} $, and $\Pi_2$ be the orthogonal projection onto the subspace spanned by $\xi_{n_1+1}, \xi_{n_1+2},\dots, \xi_n$. We set $n_2=n-n_1$.
Let $U$ be a unitary operator on $\mathcal{H}_C$, suppose $B_1=\Pi_1U$ and  $B_2=\Pi_2U$. Then \[B_{1k}=\overline{\omega}_kB_1, \;\;B_{2k}=\omega_kB_2,\] where $\omega_k=e^{ik}$ and  it follows $U_k=B_{1k}+B_{2k}$ and  \\
\begin{equation}
\mathcal{L}_{kk^{\prime}} (\rho) =p(B_{1k}\rho B^*_{1k^{\prime}}+B_{2k}\rho B^*_{2k^{\prime}})+qU_k\rho U^*_{k^{\prime}}
\end{equation} 
or
\begin{equation}
\mathcal{L}_{kk^{\prime}} (\rho) =p(\overline{\omega}_k\Pi_1U\rho U^*\Pi_1^*\omega_{k^{\prime}}+\omega_k\Pi_2U\rho U^*\Pi^*_2 \overline{\omega}_{k^{\prime}})+qU_k\rho U^*_{k^{\prime}}.
\end{equation}

In order to analyze the spectrum of the superoperator $\mathcal{L}_{kk^{\prime}}$, we will use the normalized Gell-Mann basis. Let $\; E_{jk}$ be the matrix with 1 in the $jk$-th entry and 0 elsewhere. Consider the space of $d \times d$ complex matrices, $\mathbb{C}^{d \times d}$, for a fixed $d$. Define the following matrices

* For $ k < j$, $f_{k,j} ^d = E_{kj} + E_{jk}$.

* For $ k > j$, $f_{k,j} ^d = - i ( E_{jk} - E_{kj} )$.

* For $ k = j$, Let  $f_{1,1} ^d = h_1 ^d = I_d $, the identity matrix.

* For $1 < k = j< d$, $f_{k,k} ^d =h_k ^d = h ^{d-1} _k \oplus 0$.

* For $k =j= d$, $f_{d,d} ^d =h_d ^d = \sqrt{\frac{2}{d(d-1)}} (h_1 ^{d-1} \oplus (1-d))$.\\
The collection of matrices $\{f_{k,j} ^d , \; 1\le k, j \le d\}$ are called the generalized Gell-Mann matrices in dimension $d$.  Here  $\oplus$  means Matrix Direct Sum. The generalized Gell-Mann matrices are Hermitian and traceless by construction, just like the Pauli matrices. One can also check that they are orthogonal in the Hilbert-Schmidt inner product on $\mathbb{C}^{d \times d}$. By dimension count, one sees that they span the vector space of $d\times d$  complex matrices. When $d= 2$, they are Pauli matrices. When $d=3$, they are Gell-Mann matrices.

Notes that $\{f_{k,j} ^n\}$ are not normalized. We put $\gamma_{k,j} ^n=\frac{f_{k,j} ^n}{||f_{k,j} ^n||}$, then $\{\gamma_{k,j} ^n , \; 1\le k, j \le n\}$ are orthonormal basis of $n\times n$ complex matrices in the Hilbert-Schmidt inner product.
For short notation, we set ${\gamma_l}^n=\gamma_{l,l} ^n=\frac{h_{l} ^n}{||h_{l} ^n||}$. When we fixed the dimension $n$, we omit the superscript n. 
i.e, $\gamma_{k,j}=\gamma_{k,j} ^n, \gamma_{l}=\gamma_{l}^n$. We order the basis $\gamma_{ij}$ by $\gamma_{11}, \gamma_{12}, \dots, \gamma_{1n},$ $ \gamma_{21}, \gamma_{22},\dots, \gamma_{2n}$$,\dots, \gamma_{n1}, \gamma_{n2},\dots,\gamma_{nn}.$
Now we are ready to prove


\begin{lemma}\label{le:LKK}
 Suppose $U\in U(n), \Pi_1, \Pi_2$ are defined as above and  $\{B_j\}$ is unital. Then $\mathcal{L}_{k,k+\nu}$ has the following form in terms of normalized Gell-Mann basis $\gamma_{ij}$: the first column of $\mathcal{L}_{k,k+\nu} $ satisfies 

\[
   \begin{cases}
    <\gamma_{kj}, \mathcal{L}_{k,k+\nu} (\gamma_1) >=0,\;  \forall k\neq j\\

    <\gamma_l, \mathcal{L}_{k,k+\nu} (\gamma_1) >= \begin{cases}
                                                                             \frac{2n_2\cos \nu}{n}+\frac{n_1-n_2}{n}\omega_\nu, \;l=1\\
    0, \;2\le l \le n_1\\
    \frac{2n_1i\sin \nu}{\sqrt{n(l-1)l}},  \;n_1+1 \le l \le n,
                                                                           \end{cases}

   \end{cases}
\]
where $\omega_\nu=e^{i\nu}$. The first row of $\mathcal{L}_{k,k+\nu}$ satisfies $$<\gamma_1, \mathcal{L}_{k, k+\nu}(\gamma)>=\omega_{\nu} \delta_{\gamma \gamma_1} -2i\sin \nu \frac{1}{\sqrt{n}}\;Tr({B_2}\gamma{B_2}^*),$$  where $\gamma$ is a normalized basis.
In particular, if $\nu=0$, then $\mathcal{L}_{k,k}$ has the following representation

 \[
\begin{pmatrix}
1& 0& 0& 0\\
0& \times& \times&\times\\
0& \times& \times& \times\\
0& \times& \times& \times
\end{pmatrix}.
\]
\end{lemma}
Proof: For  $n_1\le i \le n-1$,
\[\gamma_i=\frac{1}{||h_i||}
\begin{pmatrix}
1 & 0   &\cdots &0&0&\cdots &0&0 &\cdots&0 \\
0 & 1 & \cdots &0&0&\cdots &0&0&\cdots&0\\
\vdots & \vdots &\ddots & \vdots & \vdots &\ddots&\vdots&\vdots &\vdots&\vdots\\
0& 0 &\cdots&1 &0&\cdots&0 &0  &\cdots&0 \\ 
\hline
0 & 0&\cdots& 0& 1  & \dots&0 &0  &\cdots &0\\
\vdots & \vdots &\ddots & \vdots & \vdots &\ddots&\vdots&\vdots&\ddots&\vdots\\
0 & 0&\cdots& 0&0 & \cdots & 1&0 &\cdots &0 \\
0 & 0 &\cdots& 0  &0 &\cdots &0 &-(i-1)& \cdots &0\\
\vdots & \vdots &\ddots & \vdots & \vdots &\ddots&\vdots&\vdots&\ddots&\vdots\\
0 & 0 &\cdots & 0 &0&\cdots& 0&0 &\cdots &0
\end{pmatrix}\]
 where $||h_i||=\sqrt{i(i-1)}$ and there are $n_1$ rows above the line and the first $i-n_1-1$ rows below the line with one nonzero entry. It follows 
\begin{eqnarray*}
\mathcal{L}_{kk^{\prime}} (\gamma_1)&=&\frac{1}{\sqrt{n}}\left [p(\omega_{k^{\prime}-k}\Pi_1U I U^*\Pi_1^*+\omega_{k^{\prime}-k}\Pi_2U I U^*\Pi_2^*)+q(\omega_{k^{\prime}-k}\Pi_1U I U^*\Pi_2^*+\omega_{k^{\prime}-k}\Pi_2U I U^*\Pi_1^*)\right]\\
&=&\frac{1}{\sqrt{n}}(\omega_{k^{\prime}-k}\Pi_1+\omega_{k^{\prime}-k}\Pi_2)\\
&=&\frac{1}{\sqrt{n}}
\left[
\begin{array}{cccc|rrrr}
\omega_{\nu} & 0  &\cdots &0&0 &0& \cdots &0 \\
0 & \omega_\nu & \cdots &0&0 &0& \cdots &0\\
\vdots & \vdots &\vdots & \vdots & \vdots &\vdots&\vdots&\vdots\\
0& 0 & \cdots & \omega_\nu &0&0&\cdots  &0 \\ 
\hline
0 & 0& \cdots & 0& \overline{\omega}_\nu  & 0&\cdots &0\\
0 & 0 & \cdots  & 0  &0 &\overline{\omega}_\nu  & \cdots &0\\
\vdots & \vdots &\vdots & \vdots & \vdots &\vdots&\vdots&\vdots\\
0 & 0 & \cdots  & 0 &0&0& \cdots &\overline{\omega}_\nu
\end{array}\right ] 
\end{eqnarray*}
where $\nu=k^{\prime}-k$, and the block at $(1,1)$ and $(2,2)$ positions are $n_1\times n_1$ and $n_2 \times n_2$ respectively.  Hence
\begin{eqnarray*}
<\gamma_i, \mathcal{L}_{k, k+\nu} (\gamma_1) >&=&\frac{1}{\sqrt{n}}\frac{1}{\sqrt{i(i-1)}}[n_1\omega_\nu+(i-n_1-1)\overline{\omega}_\nu-(i-1)\overline{\omega}_\nu]\\
&=&\frac{1}{\sqrt{n}}\frac{1}{\sqrt{i(i-1)}}2n_1i\sin \nu.\\
\end{eqnarray*}
For 
\[\gamma_n=\frac{1}{||h_n||}
\begin{pmatrix}
1 & 0 &0   &\cdots &0&0 & \cdots &0 \\
0 & 1&0 & \cdots &0&0 & \cdots &0\\
\vdots & \vdots &\vdots & \vdots & \vdots &\vdots&\vdots&\vdots \\
0 &0&\cdots&1 &0&0&\cdots   \\ 
\hline
0 & 0& \cdots & 0& 1  & 0&\cdots &0 \\
0 & 0 & \cdots & 0  &0 &1  & \cdots &0 \\
\vdots & \vdots &\vdots & \vdots & \vdots &\vdots&\vdots\\
0 & 0 & \cdots  & 0 &0&0& \cdots &-( n-1) 
\end{pmatrix}\]
we have
\begin{eqnarray*}
<\gamma_n, \mathcal{L}_{k, k+\nu} (\gamma_1) >&=&\frac{1}{\sqrt{n}}\frac{1}{\sqrt{n(n-1)}}[n_1\omega_\nu+(n_2-1)\overline{\omega}_\nu-(n-1)\overline{\omega}_\nu]\\
&=&\frac{1}{\sqrt{n}}\frac{1}{\sqrt{n(n-1)}}2n_1i\sin \nu.\\
\end{eqnarray*}
As for $f_{kj}^n$,  the nonzero entry 1 is always off diagonal, hence
\[<\gamma_{kj}, \mathcal{L}_{k, k+\nu} (\gamma_1) >=0,\; \forall k\neq j.\] Therefore on the first column of $\mathcal{L}_{k, k+\nu} $, we have 
\[
   \begin{cases}
    <\gamma_{kj}, \mathcal{L}_{k, k+\nu} (\gamma_1) >=0,\;  \forall k\neq j\\

    <\gamma_l, \mathcal{L}_{k, k+\nu} (\gamma_1) >= \begin{cases}
                                                                             \frac{2n_2\cos \nu}{n}+\frac{n_1-n_2}{n}\omega_\nu, \;l=1\\
    0, \;2\le l \le n_1\\
    \frac{2n_1i\sin \nu}{\sqrt{n(l-1)l}},  \;n_1+1 \le l \le n.
                                                                           \end{cases}

   \end{cases}
\]
Let $\gamma$ be  a normalized basis in n-dimension, then
\small{
\begin{align*}
<\gamma_1, \mathcal{L}_{k, k+\nu}(\gamma)>
&=<\gamma_1,\;\omega_{\nu}\Pi_1U \gamma U^*\Pi_1^*+\omega_{\nu}\Pi_2U \gamma U^*\Pi_2^* >+ q<\gamma_1,\omega_{\nu}\Pi_1U \gamma U^*\Pi_2^*+\omega_{\nu}\Pi_2U \gamma U^*\Pi_1^*>\\
&=\overline{<\gamma, \overline{\omega}_{\nu}U^*{\Pi_1}^* \gamma_1{\Pi_1}U+\omega_{\nu}U^*{\Pi_2}^*\gamma_1{\Pi_2}U>}+q \overline{<\gamma, \overline{\omega}_{\nu}{\Pi_2}^* \gamma_1{\Pi_1}U+\omega_{\nu}U^*{\Pi_1}^*\gamma_1{\Pi_2}U>}\\
&=\overline{<\gamma, \overline{\omega}_{\nu}{B_1}^*{B_1} +\overline{\omega}_{\nu}{B_2}^*{B_2}+q (\overline{\omega}_{\nu}{B_2}^*B_1 +\overline{\omega}_{\nu}{B_1}^*B_2)>}\\
&=\frac{1}{\sqrt{n}} \overline{<\gamma, \overline{\omega}_{\nu} I+2i\sin \nu {B_2}^*{B_2}>}\\
&=\omega_{\nu}\frac{1}{\sqrt{n}} \overline{<\gamma,  I>} -2i\sin \nu  \frac{1}{\sqrt{n}}\overline{<\gamma, {B_2}^*{B_2}>}\\
&=\omega_{\nu} \overline{<\gamma,  \gamma_1>} -2i\sin \nu  \frac{1}{\sqrt{n}}<{B_2}^*{B_2}, \gamma>\\
&=\omega_{\nu} \delta_{\gamma \gamma_1} -2i\sin \nu  \frac{1}{\sqrt{n}}<{B_2}^*{B_2}, \gamma>.
\end{align*}}
Therefore on the first row of $\mathcal{L}_{k, k+\nu}$ in Gell-Mann basis:\\
$$<\gamma_1, \mathcal{L}_{k, k+\nu}(\gamma)>=\omega_{\nu} \delta_{\gamma \gamma_1} -2i\sin \nu \frac{1}{\sqrt{n}}\;Tr({B_2}\gamma{B_2}^*).$$
By plug in $\gamma=\gamma_i, i=1,2,\dots, n$ and $\nu=0$, the matrix representation of $\mathcal{L}_{k,k+\nu}$ has the desired form listed in the lemma. In particular, if $h=\gamma_1$, then \[<{B_2}^*{B_2}, \gamma_1>= Tr (\Pi_2 U \gamma_1 U^*{\Pi_2}^*)=\frac{1}{\sqrt{n}} Tr (\Pi_2{\Pi_2}^*)=\frac{n_2}{\sqrt{n}},\] we have 
\[l_{11}(\nu)= \omega_{\nu} -2i\sin \nu \frac{n_2}{n}= \omega_{\nu} \frac{n_1}{n}+\overline{\omega}_{\nu} \frac{n_2}{n},\]and \[l_{11}(0)=1.\]


\begin{theorem}\label{CONV}
 Suppose $U\in U(n), \Pi_1, \Pi_2$ are defined as above. If $1$ is an eigenvalue of $\mathcal{L}_{kk}$ with algebraic multiplicity 1 and $|\lambda|<1 $ for any other eigenvalue $\lambda$ of $\mathcal{L}_{kk}$.Then

\begin{enumerate}
\item[(a)] $z^{\prime}_0(0)=\frac{(n_2-n_1)i}{n}$\\
\item[(b)] $\lim_{t\rightarrow\infty}\hat{P}(\frac{\nu}{t},t)=e^{-z_0^{\prime}(0)\nu}$\\
\item[(c)] $ \lim_{t\rightarrow\infty}\hat{P}(\frac{\nu}{\sqrt{t}},t)e^{z_0^{\prime}(0)\nu \sqrt{t}}= \frac{1}{2\pi}\int_0^{2\pi}  e^{-\frac{1}{2}[z^{\prime\prime}_0(0)-(z_0^{\prime}(0))^2]\nu^2} dk$\\
\end{enumerate}
where $z_0(\nu)$ is the root of $det (1-z\mathcal{L}_{k,k+\nu})=0$ such that $z_0(0)=1$.
\end{theorem}
Proof. a) We order the generalized Gell-Mann orthonormal basis $\gamma_{ij}$ by $\gamma_{11}, \gamma_{12}, \dots, \gamma_{1n},$ $ \gamma_{21}, \gamma_{22},\dots, \gamma_{2n}$$,\dots, \gamma_{n1}, \gamma_{n2},\dots,\gamma_{nn}$. Let $A$ be the $m\times m$ matrix representation of $1-z\mathcal{L}_{k,k}$ with respect to $\{\gamma_{ij}\}$ in the given order. 
Any operator $\hat{O}$ which act on $\mathcal{H}_C$ can be represented by linear combination of  matrices $\{\gamma_{ij}\}$:
\begin{equation}
\hat{O}=s_1\gamma_{11}+s_2 \gamma_{12}+\cdots+s_{m}\gamma_{nn}, \quad m=n^2.
\end{equation}
Hence  $\hat{O}$ can be represented by a column vector $\hat{O}=(s_1, s_2,\dots, s_{m})^T$. Let $A_{ij}$ be the cofactor of $A$ at the $ij$-th entry, then
\begin{equation*}
A^{-1}\hat{O}=\frac{1}{\det (A)}\begin{pmatrix}
A_{11} & A_{21} &\cdots &A_{m1}\\
A_{12} & A_{22} &\cdots &A_{m2}\\
\vdots & \vdots &\cdots&\vdots\\
A_{1m} & A_{2m} &\cdots &A_{mm}
\end{pmatrix}\begin{pmatrix}
s_1\\
s_2\\
\vdots\\
s_{m}
\end{pmatrix}.
\end{equation*}
 Since $Tr (\gamma_{ij})=0,$ for all $i,j$ except for $Tr(\gamma_{11})=\sqrt{n},$ we have 
\begin{equation}
\hat{P}(\nu,t)<e^{i\nu x}>_t=\frac{1}{2\pi i}\oint_{|z|=r<1} \frac{G(z,v)}{z^{t+1}}dz=\frac{1}{2\pi}\int dk\frac{1}{2\pi i}\oint_{|z|=r<1}\frac{\sqrt{n}h(z,\nu)} {z^{t+1} det A },
\end{equation}
   for some $0<r<1$, where 
   $h(z,v)=A_{11}s_1+A_{21}s_2+\cdots+A_{m1}s_{m}$ using the same argument as in  [\cite{FFXY}, Theorem (3.1)], we have 

\[
\lim_{t\rightarrow\infty}\hat{P}(\frac{\nu}{t},t)=\frac{\sqrt{n}h(1,0)}{\frac{\partial g}{\partial z}(1,0)}\lim_{t\rightarrow\infty}z_0\left(\frac{\nu}{t}\right)^{-t-1},
\]
where $g(z,\nu)=det(A)$. If $\nu=0,$ the first row of $A$ is equal to $0$ except $a_{11}=1-z$. We will show that $\frac{\sqrt{n}h(1,0)}{\frac{\partial g}{\partial z}(1,0)}=-1$  and $z^{\prime}_0(0)=\frac{(n_1-n_2)i}{n}$. Note that $z^{\prime}_0(0)=\frac{(n_1-n_2)i}{n}$ is independent of $k$. 

Let's first  show that $\frac{\sqrt{n}h(1,0)}{\frac{\partial g}{\partial z}(1,0)}=-1$. Let $M$ denote the following submatrix of $L=(l_{ij})$:
\[
M(\nu)=\begin{pmatrix}
l_{22}&\cdots&l_{2m}\\
l_{32}&\cdots&l_{3m}\\
\vdots&\vdots&\vdots\\
l_{m2}&\cdots&l_{mm}
\end{pmatrix}.
\]
Then the matrix \[A|_{\nu=0}=
\begin{pmatrix}
1-z&\textbf{0} \\
\textbf{0}&I_{m-1}-zM(0)
\end{pmatrix}.
\]
The cofactors $A_{21}=A_{31}=\cdots=A_{m1}=0$
and $A_{11}=\det(I_{m-1}-zM(0))$.
By Lemma \ref{le:LKK},
 $\frac{1}{z_i(0)}$ are eigenvalues of $\mathcal{L}_{k,k}$ for $i=0,1,\dots, m-1$. Therefore $\frac{1}{z_1(0)},\frac{1}{z_2(0)},\dots,\frac{1}{z_{m-1}(0)}$ are eigenvalues of $M(0)$.
   Hence
\begin{eqnarray*}
\det(I_{m-1}-zM(0))&=&\left(1-\frac{z}{z_1(0)}\right)\left(1-\frac{z}{z_2(0)}\right)\left(1-\frac{z}{z_{m-1}(0)}\right)\\
&=&-\frac{(z-z_1(0))(z-z_2(0))\cdots(z-z_{m-1}(0))}{z_1(0)z_2(0)\cdots z_{m-1}(0)}.
\end{eqnarray*}
On the other hand,
\begin{eqnarray}
\frac{\partial g}{\partial z}(1,0)=\frac{(1-z_1(0))(1-z_2(0))\cdots (1-z_{m-1}(0))}{z_1(0)z_2(0)\cdots z_{m-1}(0)}.\label{eq:partial gz}
\end{eqnarray}
Hence
\begin{eqnarray*}
\frac{\sqrt{n}h(1,0)}{\frac{\partial g}{\partial z}(1,0)}=\sqrt{n}s_1\frac{-\frac{1}{z_1(0)z_2(0)\cdots z_{m-1}(0)}(1-z_1(0))(1-z_2(0))\cdots(1-z_{m-1}(0))}{\frac{1}{z_1(0)z_2(0)\cdots z_{m-1}(0)}(1-z_1(0))(1-z_2(0))\cdots (1-z_{m-1}(0))}\\
=-\sqrt{n}s_{1}=-\sqrt{n}<\gamma_{11},\hat{O}>=-\sqrt{n}Tr(\frac{1}{\sqrt{n}}I \hat{O})=-Tr(\hat{O})=-1,
\end{eqnarray*}
since $\hat{O}$ is a density operator.

Next we  show that $z_0^{\prime}(0)=\frac{(n_1-n_2)i}{n}$. Since $g(z_0(\nu),\nu)=0$, we have
\begin{eqnarray}\label{generating}
0=\frac{d}{d\nu} g(z_0(\nu),\nu)=\frac{\partial}{\partial \nu} g(z,\nu)\Big|_{z=z_0(\nu)}+\frac{\partial}{\partial z} g(z,\nu)\Big|_{z=z_0(\nu)}z_0^{\prime}(\nu).
\end{eqnarray}
When  $\nu=0$, we have $z_0(0)=1$ and it follows that 
\begin{equation}
0=\frac{\partial}{\partial \nu} g(1,\nu)\Big|_{\nu=0}+\frac{\partial}{\partial z}\ g(z,0)\Big|_{z=1}z_0^{\prime}(0). \label{eq:implicity derivative}
\end{equation}
By the form of the first column of $\mathcal{L}_{k,k}$, we have 
\[g(1,\nu)=\left(1-  \frac{2n_2\cos \nu}{n}-\frac{n_1-n_2}{n}\omega_\nu \right)A_{11}\nu
-\sum_{l=n_1+1}^{n}   \frac{2n_1i\sin \nu}{\sqrt{n(l-1)l}} A_{l1}\nu.\]
By Lemma \ref{le:LKK}, the cofactor  $A_{l1}(0)=0, l=n_1+1,\dots, n$. It follows that
\small{
\begin{eqnarray*}
\frac{\partial}{\partial \nu} g(1,\nu)\Big|_{\nu =0}&=&\left(\frac{2n_2\sin \nu}{n}-\frac{n_1-n_2}{n}i\omega_\nu \right)A_{11}\nu +\left(1-  \frac{2n_2\cos \nu}{n}-\frac{n_1-n_2}{n}\omega_\nu \right){A^{\prime}}_{11}\nu\\
&-&\sum_{l=n-1+1}^{n}   \frac{2n_1i\cos \nu}{\sqrt{n(l-1)l}} A_{l1}\nu-\sum_{l=n-1+1}^{n}   \frac{2n_1i\sin \nu}{\sqrt{n(l-1)l}} {A^{\prime}}_{l1}\nu\\ \label{eq:partial gnu}\\
&=& \frac{-(n_1-n_2)i}{n}A_{11}(0)=\frac{-(n_1-n_2)i}{n}det[1-M(0)].
\end{eqnarray*}
}On the other hand, by
\[\frac{\partial g(z,0)}{\partial z}\Big|_{z =1}=-det[1-M(0)], \]
it follows that 
$z^{\prime}_0(0)=\frac{-(n_1-n_2)i}{n}.$  Therefore 
\\
\[\lim_{t\rightarrow\infty}\hat{P}\left(\frac{\nu}{t},t\right)= \frac{1}{2\pi}\int_0^{2\pi}dk \lim_{t\rightarrow\infty}[1+z_0^{\prime}(0)\frac{\nu}{t} +o(\frac{\nu}{t})] ^{-t-1}=e^{-z_0^{\prime}(0)\nu}=e^{\frac{(n_1-n_2)i\nu}{n}}.\]
 Let $\phi(\nu)=z_0\left(\nu \right) \exp\left({-z^{\prime}(0)\nu }\right)$, then $\phi(0)=1, \phi^{\prime}(0)=0$ and $\phi^{\prime\prime}(0)=z_0^{\prime\prime}(0)-[z_0^{\prime}(0)]^2$, hence
\begin{align*}
 \lim_{t\rightarrow\infty}\hat{P}\left(\frac{\nu}{\sqrt{t}},t\right)\exp\left({z^{\prime}(0)\nu \sqrt{t}}\right)&=\lim_{t\rightarrow\infty}z_0\left(\frac{\nu}{\sqrt{t}}\right)^{-t-1}\exp\left({z^{\prime}(0)\nu \sqrt{t}}\right)\\
&=\lim_{t\rightarrow\infty}\left [z_0\left(\frac{\nu}{\sqrt{t}}\right) \exp\left({-z^{\prime}(0)\frac{\nu}{\sqrt{t}}}\right)\right]^{-t} \lim_{t\rightarrow\infty}z_0^{-1}\left(\frac{\nu}{\sqrt{t}}\right)\\
&= \frac{1}{2\pi}\int_0^{2\pi}dk \lim_{t\rightarrow\infty}[1+\phi_0^{\prime\prime}(0)\frac{\nu^2}{t} +o(\frac{\nu}{t})] ^{-t}\\
&= \frac{1}{2\pi}\int_0^{2\pi}  e^{-\frac{1}{2}[z^{\prime\prime}_0(0)-(z_0^{\prime}(0))^2]\nu^2} dk.
\end{align*}

\section{Superoperators and the eigenvalue condition}
\setcounter{equation}{0}

In this section, we prove a Perron-Frobenius type of theorem which is very useful for determining if the superoperator $\mathcal{L}_{kk}$ satisfies the eigenvalue condition. That is, one is its only eigenvalue with absolute value 1 and the absolute values of all other eigenvalues are strictly less than 1.\\ 
\begin{theorem} \label{EigenCond}
 Let $\mathcal{H}$ be a finite dimensional Hilbert space over the complex field $\mathbb{C}$ equipped with the regular inner product, and  $A$ and $B$ are linear operators on $\mathcal{H}$.  Let $\mathcal{L}=qA+pB, 0<p,q<1, p+q=1$. Suppose i) $||A||, ||B||\le 1;$ ii) 1 is an eigenvalue of A and B and there is a $\rho_1$ such that $A^*\rho_1=\rho_1, B^*\rho_1=\rho_1$, then\\ 
 a) If $\mathcal{L}\rho=\lambda\rho$ with $|\lambda|=1$, then $A\rho=\lambda\rho$ and $B\rho=\lambda\rho$.\\
b) If the multiplicity of 1as an eigenvalue of B is 1, then the multiplicity of 1 as an eigenvalue of $\mathcal{L}$ is also 1.

\end{theorem}

\begin{corollary}
Under the hypothesis of Theorem (\ref{EigenCond}). If $\lambda=1$ is the only eigenvalue of B with $|\lambda|=1$, then $\lambda=1$ is also the only eigenvalue of $\mathcal{L}$ with $|\lambda|=1$.

\end{corollary}


Proof of a):  $\mathcal{L}\rho=\lambda\rho$ with $|\lambda|=1$ implies that 
\[||\rho||=|\lambda|||\rho||=||qA\rho+pB\rho|| \le ||qA\rho||+||pB\rho||=q||A||||\rho||+p||B||||\rho||\le q||\rho||+p||\rho||=||\rho||.\]
The equality holds only when $A\rho=\alpha B\rho$ for some $\alpha>0$.
The assumptions i) and ii) imply that $||A||=||B||=1$, since $\rho \ne 0$, so it follows that $\alpha=1$ and  $A\rho=\lambda\rho$ and $B\rho=\lambda\rho.$\\

Proof of b)  Suppose $\rho$ is a generalized eigenvector of $\mathcal{L}$, then there is an $l\ge 2$ such that $(\mathcal{L}-1)^l\rho=0$ or $(\mathcal{L}-1)(\mathcal{L}-1)^{l-1}\rho=0$. It follows that $(\mathcal{L}-1)^{l-1}\rho$ is  an eigenvector of $\mathcal{L}$ with eigenvalue 1. By a), $(\mathcal{L}-1)^{l-1}\rho$ is also an eigenvector of $B$ with eigenvalue 1.  Since the multiplicity of 1 as an eigenvalue of B is 1,  the dimension of the eigenspace of 1 is 1, and we have $(\mathcal{L}-1)^{l-1}\rho=\beta\rho_0$ for some $\beta\in \mathbb{C}$, where $\rho_0$ is an eigenvector of B with eigenvalue  1.   We will show that $(\mathcal{L}-1)^{l-1}\rho=0$ for the next step. Note that
\begin{eqnarray*}
 \beta<\rho_1, \rho_0>=<\rho_1, \beta\rho_0>&=&<\rho_1, (\mathcal{L}-1)^{l-1}\rho>=<\rho_1, (\mathcal{L}-1)(\mathcal{L}-1)^{l-2}\rho>\\
 &=&<\rho_1, \mathcal{L} (\mathcal{L}-1)^{l-2}\rho>-<\rho_1, (\mathcal{L}-1)^{l-2}\rho>\\
 &=&<\mathcal{L} ^*\rho_1, (\mathcal{L}-1)^{l-2}\rho>-<\rho_1, (\mathcal{L}-1)^{l-2}\rho>\\
 &=&<\rho_1, (\mathcal{L}-1)^{l-2}\rho>-<\rho_1, (\mathcal{L}-1)^{l-2}\rho>\\
 &=&0,
 \end{eqnarray*}
it is sufficient to show that $\beta=0$ or equivalently $<\rho_1, \rho_0>\ne 0$.
 Since the dimension of the eigenspace of 1 as an eigenvalue of B is 1, it is also true for $B^*$. If   $<\rho_1, \rho_0>= 0, $ then $\rho_0 \perp \text{Ker}(B^*-I)$ or $\rho_0\in \overline{\text{Ran}(B-I)}=\text{Ran}(B-I)$. The last equality holds since B is a bounded linear and therefore continuous operator on $\mathcal{H}$. Therefore $\rho_0=(B-I)x$ for some $x\in \mathcal{H}$ and  $(B-I)^2x=(B-I)\rho_0=0$, that is, $x$ is a generalized eigenvector of B of eigenvalue 1. Then $x=\gamma \rho_0$ for some $\gamma \in \mathbb{C}$ and we have  $\rho_0=(B-I)x=\gamma (B-I)\rho=0$, which is a contradiction. Therefore $<\rho_1, \rho_0>\ne 0$ and $(\mathcal{L}-1)^{l-1}\rho=0$.\smallskip
 
 Repeat the above argument,  we have $(\mathcal{L}-I)\rho=0$, that is, $\rho$ is a genuine eigenvector of $\mathcal{L}$ of eigenvalue 1. The multiplicity of 1 as an eigenvalue of $\mathcal{L}$ is 1. \\
 
 Proof of Corollary:  the assumption ii) implies that $\mathcal{L}^*\rho_1=\rho_1$ or \[1\in Spec(\mathcal{L}^*)=Spec(\mathcal{L}).\]On the other hand, if $\lambda \in Spec(\mathcal{L})$ with $|\lambda|=1$ and $\lambda \ne 1$, then by Theorem \ref{EigenCond}(a), $\lambda$ is also an eigenvalue of B,  a contradiction. 
 
 
 The following theorem is a special version of Theorem (\ref{EigenCond}) which applies to the partially open quantum random walk on $Z^1$ defined in this paper.

 \begin{theorem}\label {eigencond}
 Let $\mathcal{H}$ be a finite dimensional Hilbert space over the complex field $\mathbb{C}$ equipped with the regular inner product. $C(\rho)=U\rho U^*$ and $D(\rho)=\sum_{j=1}^{n}A_j\rho{A_j}^*$ are linear operators on $\mathcal{H}$ for $\rho \in \mathcal{H}$. Where $U$ is a unitary operator and $\{ A_j\}$ are unital measurements. Define $\mathcal{L}(\rho)=qC(\rho)+pD(\rho), 0<p,q<1, p+q=1$. Suppose
$\lambda=1$ is the only eigenvalue of $D(\rho)$ with $|\lambda|=1,$ 
then  $\lambda=1$ is the only eigenvalue of $\mathcal{L}$with $|\lambda|=1$.  Moreover, if the multiplicity of $\lambda=1 $ as an eigenvalue of D is 1, then the multiplicity of $\lambda=1$ as an eigenvalue of $\mathcal{L}$ is also 1. 
\end{theorem}
 
 Proof:  Let $\rho_1=\frac{I_{\mathcal{H}}}{dim(\mathcal{H})}$, then it is easy to verify  that $C(\rho_1)=U\rho_1 U^*$ and $D(\rho_1)=\sum_{j=1}^{n}A_j\rho_1{A_j}^*$ satisfy the assumption of Theorem \ref{EigenCond} and therefore   theorem  follows.

 \begin{theorem}\label{eigencond1}
Consider a decoherent quantum random walk on $Z^1$ with coin space $\mathcal{H}_C$  given by $\mathcal{L}_{k,k^{\prime},p}(\rho)=qC(\rho)+pD(\rho), 0<p\le 1, p+q=1$. Where  $C(\rho)=U_k\rho U^*_{k^{\prime}}$ and $D(\rho)=\sum_{j=1}^{m}\Pi_j U_k\rho{(\Pi_j}U_{k^{\prime}})^*$ are linear operators on $\mathcal{H}$ for $\rho \in L(\mathcal{H})$ and $U_k$ is a unitary operator. Assume
 $\mathcal{L}_{k,k,1}$ satisfies the eigenvalue condition, then $\mathcal{L}_{k,k, p}$ satisfies the eigenvalue condition for any $1\ge p>0$.
\end{theorem}

\begin{remark}
 $\mathcal{L}_{k,k^{\prime},0}$  is the coherent quantum random walk and 
$\mathcal{L}_{k,k^{\prime},1}$  is the decoherent quantum random walk.
\end{remark}
 
 Proof:   Assume $\mathcal{L}_{k,k, 1}$ satisfies eigenvalue condition.  By Theorem (\ref{eigencond}), $\mathcal{L}_{k,k, p}$ satisfies eigenvalue condition for all $0<p<1$ and therefore for all $0<p\le 1$. \\


\section{Applications}

In this section, we will illustrate two examples for $n=2, 3$ and compute the limiting distributions of characteristic equations of the position probability functions $p(x,t).$
First we find a general formula to compute the variance $\sigma^2=z^{\prime\prime}_0(0)-z^{\prime}(0)^2$. Note that $g(z, \nu)=0$, direct  computation shows that
\begin{align*}
0=&\frac{d }{d\nu}\Big[\frac{\partial g(z,\nu)}{\partial \nu}\Big\vert _{z=z_0(\nu)}+\frac{\partial g(z,\nu)}{\partial z}\Big\vert_{z=z_0(\nu)}z_0^{\prime}(\nu)\Big]\\
=&\frac{\partial g_{\nu}(z,\nu)}{\partial \nu}\Big\vert_{z=z_0(\nu)}+\frac{\partial g_{\nu}(z,\nu)}{\partial z}\Big\vert_{z=z_0(\nu)}z_0^{\prime}(\nu) \\
&+
\Big(\frac{\partial g_z(z,\nu)}{\partial \nu}\Big\vert_{z=z_0(\nu)}+\frac{\partial g_z(z,\nu)}{\partial z}\Big\vert_{z=z_0(\nu)}z_0^{\prime}(\nu)\Big)z_0^{\prime}(\nu)  +g_z(z_0(\nu),\nu)z_0^{\prime\prime}(\nu).
\end{align*}
When $\nu=0$, we have $z_0(\nu)=1$. It  then follows that
\[g_{\nu\nu}(1,0)+g_{\nu z}(1,0)z_0^{\prime}(0)+g_{z\nu}(1,0)z_0^{\prime}(0)+g_{zz}(1,0)(z_0^{\prime}(0))^2+g_{z}(1,0)z_0^{\prime\prime}(0)=0.\]
Hence
\[\sigma^2=z^{\prime\prime}_0(0)-z^{\prime}(0)^2=\Big[-\frac{g_{zz}(1,0)}{g_z(1,0)}-1\Big]z_0^{\prime}(0)^2-\frac{g_{z\nu}(1,0)}{g_z(1,0)}z_0^{\prime}(0)-\frac{g_{\nu\nu}(1,0)}{g_z(1,0)}.\]
Recall  $g(z,0)=(1-z)det(I-zM(0))$. It then follows that 
\begin{align*}
g_z(z,0)&=-det(I-zM(0))+(1-z)[det(I-zM(0))]^{\prime},\\
g_{zz}(z,0)&=-2[det(I-zM(0))]^{\prime}+(1-z)[det(I-zM(0))]^{\prime\prime}.
\end{align*}
So  \;$g_z(1,0)=-det(I-M(0))$ and  $g_{zz}(1,0)=-2[det(I-M(0))]^{\prime}$.

\begin{example}

Let's consider  the general  Hadamard walk on $Z^1$,  the evolution operator  is
\begin{equation}
U_k=
\begin{pmatrix}
\overline{\omega}_{k} \cos \theta&\overline{\omega}_{k} \sin \theta\\
\omega_{k} \sin \theta &-\omega_{k} \cos \theta
\end{pmatrix}.
\end{equation}
Suppose $\Pi_1U_k=\overline{\omega}_{k} 
\begin{pmatrix}
 \cos \theta&\sin \theta\\
0 &0
\end{pmatrix}$ 
and 
$\Pi_2U_k=\omega_{k} 
\begin{pmatrix}
0 &0\\
\sin \theta&-\cos \theta\\

\end{pmatrix}$.  Order the Pauli matrices as in the generalized Gell-Mann matrices:\\
$\gamma_{11}=\frac{1}{\sqrt{2}}
\begin{pmatrix}
1&0\\
0&1
\end{pmatrix}
$,\quad $\gamma_{12}=\frac{1}{\sqrt{2}}
\begin{pmatrix}
0&1\\
1&0
\end{pmatrix}
$,\quad$\gamma_{21}=\frac{-i}{\sqrt{2}}
\begin{pmatrix}
0&1\\
-1&0
\end{pmatrix}
$, \quad $\gamma_{22}=\frac{1}{\sqrt{2}}
\begin{pmatrix}
1&0\\
0&-1
\end{pmatrix}.
$
Under the Pauli matrices basis:
\[
\mathcal{L}_{k,k+\nu}=
\begin{pmatrix}
\cos \nu &i\sin\nu\sin2 \theta & 0 &i\sin\nu \cos 2\theta\\
0& -q\cos (2k+\nu)\cos 2\theta& q\sin(2k+\nu)&q\cos(2k+\nu)\sin 2\theta\\
0&-q\sin (2k+\nu)\cos2\theta&-q\cos(2k+\nu)&q\sin(2k+\nu)\sin 2\theta\\
i\sin\nu&\cos\nu\sin 2\theta&0&\cos\nu\cos2\theta
\end{pmatrix}\]
The superoperator $\mathcal{L}_{k,k+\nu}$ satisfies the eigenvalue condition (see c.f. Section 4 in \cite{FFXY}). 
 We have \[z_0^{\prime}(0)=0;\;\;
z_0^{\prime\prime}(0)=\frac{1+q\cos(2k)+q^2}{1-q^2}\cot^2\theta,
\]
where $q=1-p$, and
\[
\hat{P}(\frac{\nu}{\sqrt{t}},t)\rightarrow\frac{1}{2\pi}\int_0^{2\pi}e^{-\frac{1}{2}\frac{1+q\cos(2k)+q^2}{1-q^2}\cot^2\theta\nu^2}dk.
\]
\end{example}
\begin{example}
Let  $r=\frac{1}{\sqrt{2}},$
$U=r
\begin{pmatrix}
1 & 0 &1\\
r &1 & -r\\
-r &1& r
\end{pmatrix}$, then
the evolution operator  is given by 
\[U_k=r
\begin{pmatrix}
 \overline{\omega}_{k}  & 0 & \overline{\omega}_{k} \\
r \overline{\omega}_{k} & \overline{\omega}_{k} & -r \overline{\omega}_{k} \\
-r  \omega_{k} & \omega_{k}& r \omega_{k}
\end{pmatrix}.\]
Define 
\[B_{1k}=\Pi_1U_k=r
\begin{pmatrix}
 \overline{\omega}_{k}  & 0 & \overline{\omega}_{k} \\
r \overline{\omega}_{k} & \overline{\omega}_{k} & -r \overline{\omega}_{k} \\
0 &0&0
\end{pmatrix},\]
\[B_{2k}=\Pi_2U_k=r
\begin{pmatrix}
0 &0&0\\
0 &0&0\\
-r  \omega_{k} & \omega_{k}& r \omega_{k}
\end{pmatrix}.\]

Under the ordered Gell-Mann basis,
$\mathcal{L}_{k, k+\nu}=$
\[\tiny{
\begin{pmatrix}
\cos \nu +\frac{1}{3}i\sin \nu& \frac{1}{\sqrt{3}}i\sin \nu &0   &\frac{1}{2\sqrt{6}}i\sin \nu &-\frac{1}{\sqrt{6}}i\sin \nu&0 &-\frac{2}{\sqrt{3}}i\sin \nu &0 & -\frac{1}{2\sqrt{3}}i\sin \nu\\
0 & \frac{1}{2}\omega_{\nu} & 0 & 0 & \frac{1}{2\sqrt{2}}\omega_{\nu} &\frac{1}{2}\omega_{\nu} &  0 & 0 & \frac{3}{2\sqrt{6}}\omega_{\nu} \\
0 &\frac{q\cos{(2k+\nu)}}{2} & 0& -\frac{q\sin{(2k+\nu)}}{2} &-\frac{q\cos{(2k+\nu)}  }{2\sqrt{2}} &\frac{q\cos{(2k+\nu)}  }{2}& -\frac{q\sin{(2k+\nu)}}{\sqrt{2}} &\frac{q\sin(2k+\nu)}{2} &-\frac{\sqrt{6}q\cos{(2k+\nu)}}{4} \\
0 & 0 &0 &-\frac{1}{2}\omega_{\nu}  &0&0 &\frac{1}{\sqrt{2}}\omega_{\nu} &\frac{1}{2}\omega_{\nu} &0 \\
0 &-\frac{1}{2\sqrt{2}}\omega_{\nu} & \frac{3}{4}\omega_{\nu} & 0 & \frac{3}{8}\omega_{\nu} & \frac{1}{\sqrt{2}}\omega_{\nu}   & 0 & 0& \frac{-3}{4\sqrt{12}}\omega_{\nu}  \\
0 &0 &\frac{q\cos{(2k+\nu)}}{2}  & -\frac{q\sin{(2k+\nu)}}{\sqrt{2}} &  -\frac{3q\cos{(2k+\nu)}}{4} & 0 &  0 & -\frac{q\sin{(2k+\nu)}}{\sqrt{2}} &  \frac{3q\cos{(2k+\nu)}}{2\sqrt{12}} \\
0 &\frac{-qi\overline{\omega}_{2k+\nu}}{4}& 0&-\frac{q\cos{(2k+\nu)}}{2} & \frac{q\sin{(2k+\nu)}  }{2\sqrt{2}} &-\frac{q\sin{(2k+\nu)}  }{2} &-\frac{q\cos{(2k+\nu)}  }{\sqrt{2}}  & \frac{q\cos{(2k+\nu)}}{2}& \frac{\sqrt{6}q\sin{(2k+\nu)}}{4}  \\
0 &0 & -\frac{q\sin(2k+\nu)}{2}&-\frac{q\cos{(2k+\nu)}  }{\sqrt{2}}  &\frac{3q\sin{(2k+\nu)}  }{4}  &0&  0& -\frac{q\cos{(2k+\nu)}  }{\sqrt{2}} &-\frac{3q\sin{(2k+\nu)}  }{4\sqrt{12}}  \\
 \frac{2\sqrt{2}}{3}i\sin \nu& \frac{\omega_{\nu}}{2\sqrt{6}}+ \frac{\overline{\omega}_{\nu}}{\sqrt{6}} &  \frac{\omega_{\nu}}{4\sqrt{3}}+ \frac{\overline{\omega}_{\nu}}{2\sqrt{3}}  & 0 & \frac{\omega_{\nu}}{4\sqrt{12}}+ \frac{\overline{\omega}_{\nu}}{2\sqrt{12}} &  \frac{-\omega_{\nu}}{\sqrt{6}}- \frac{2\overline{\omega}_{\nu}}{\sqrt{6}}& 0 & 0& \frac{-\omega_{\nu}}{24}+ \frac{\overline{\omega}_{\nu}}{4}  \\
\end{pmatrix}}\]

Denote 

\[I(\lambda, \nu)=\tiny{
\begin{pmatrix}

 \frac{1}{2}\omega_{\nu} & 0 & 0 & \frac{1}{2\sqrt{2}}\omega_{\nu} &\frac{1}{2}\omega_{\nu} &  0 & 0 & \frac{3}{2\sqrt{6}}\omega_{\nu} \\
\frac{q\cos{(2k+\nu)}}{2} & 0& -\frac{q\sin{(2k+\nu)}}{2} &-\frac{q\cos{(2k+\nu)}  }{2\sqrt{2}} &\frac{q\cos{(2k+\nu)}  }{2}& -\frac{q\sin{(2k+\nu)}}{\sqrt{2}} &\frac{q\sin(2k+\nu)}{2} &-\frac{\sqrt{6}q\cos{(2k+\nu)}}{4} \\
0& 0 &-\frac{1}{2}\omega_{\nu}  &0&0 &\frac{1}{\sqrt{2}}\omega_{\nu} &\frac{1}{2}\omega_{\nu} &0 \\
-\frac{1}{2\sqrt{2}}\omega_{\nu} & \frac{3}{4}\omega_{\nu} & 0 & \frac{3}{8}\omega_{\nu} & \frac{1}{\sqrt{2}}\omega_{\nu}   & 0 & 0& \frac{-3}{4\sqrt{12}}\omega_{\nu}  \\
0&\frac{q\cos{(2k+\nu)}}{2}  & -\frac{q\sin{(2k+\nu)}}{\sqrt{2}} &  -\frac{3q\cos{(2k+\nu)}}{4} & 0 &  0 & -\frac{q\sin{(2k+\nu)}}{\sqrt{2}} &  \frac{3q\cos{(2k+\nu)}}{2\sqrt{12}} \\
\frac{-qi\overline{\omega}_{2k+\nu}}{4}& 0&-\frac{q\cos{(2k+\nu)}}{2} & \frac{q\sin{(2k+\nu)}  }{2\sqrt{2}} &-\frac{q\sin{(2k+\nu)}  }{2} &-\frac{q\cos{(2k+\nu)}  }{\sqrt{2}}  & \frac{q\cos{(2k+\nu)}}{2}& \frac{\sqrt{6}q\sin{(2k+\nu)}}{4}  \\
0 & -\frac{q\sin(2k+\nu)}{2}&-\frac{q\cos{(2k+\nu)}  }{\sqrt{2}}  &\frac{3q\sin{(2k+\nu)}  }{4}  &0&  0& -\frac{q\cos{(2k+\nu)}  }{\sqrt{2}} &-\frac{3q\sin{(2k+\nu)}  }{4\sqrt{12}}  \\
 \frac{\omega_{\nu}}{2\sqrt{6}}+ \frac{\overline{\omega}_{\nu}}{\sqrt{6}} &  \frac{\omega_{\nu}}{4\sqrt{3}}+ \frac{\overline{\omega}_{\nu}}{2\sqrt{3}}  & 0 & \frac{\omega_{\nu}}{4\sqrt{12}}+ \frac{\overline{\omega}_{\nu}}{2\sqrt{12}} &  \frac{-\omega_{\nu}}{\sqrt{6}}- \frac{2\overline{\omega}_{\nu}}{\sqrt{6}}& 0 & 0& \frac{-\omega_{\nu}}{24}+ \frac{\overline{\omega}_{\nu}}{4}  \\
\end{pmatrix}},\]

\[J(\lambda, \nu)=\tiny{
\begin{pmatrix}
 \frac{1}{\sqrt{3}}i\sin \nu &0   &\frac{1}{2\sqrt{6}}i\sin \nu &-\frac{1}{\sqrt{6}}i\sin \nu&0 &-\frac{2}{\sqrt{3}}i\sin \nu &0 & -\frac{1}{2\sqrt{3}}i\sin \nu\\
 \frac{1}{2}\omega_{\nu} & 0 & 0 & \frac{1}{2\sqrt{2}}\omega_{\nu} &\frac{1}{2}\omega_{\nu} &  0 & 0 & \frac{3}{2\sqrt{6}}\omega_{\nu} \\
\frac{q\cos{(2k+\nu)}}{2} & 0& -\frac{q\sin{(2k+\nu)}}{2} &-\frac{q\cos{(2k+\nu)}  }{2\sqrt{2}} &\frac{q\cos{(2k+\nu)}  }{2}& -\frac{q\sin{(2k+\nu)}}{\sqrt{2}} &\frac{q\sin(2k+\nu)}{2} &-\frac{\sqrt{6}q\cos{(2k+\nu)}}{4} \\
 0 &0 &-\frac{1}{2}\omega_{\nu}  &0&0 &\frac{1}{\sqrt{2}}\omega_{\nu} &\frac{1}{2}\omega_{\nu} &0 \\
-\frac{1}{2\sqrt{2}}\omega_{\nu} & \frac{3}{4}\omega_{\nu} & 0 & \frac{3}{8}\omega_{\nu} & \frac{1}{\sqrt{2}}\omega_{\nu}   & 0 & 0& \frac{-3}{4\sqrt{12}}\omega_{\nu}  \\
0 &\frac{q\cos{(2k+\nu)}}{2}  & -\frac{q\sin{(2k+\nu)}}{\sqrt{2}} &  -\frac{3q\cos{(2k+\nu)}}{4} & 0 &  0 & -\frac{q\sin{(2k+\nu)}}{\sqrt{2}} &  \frac{3q\cos{(2k+\nu)}}{2\sqrt{12}} \\
\frac{-qi\overline{\omega}_{2k+\nu}}{4}& 0&-\frac{q\cos{(2k+\nu)}}{2} & \frac{q\sin{(2k+\nu)}  }{2\sqrt{2}} &-\frac{q\sin{(2k+\nu)}  }{2} &-\frac{q\cos{(2k+\nu)}  }{\sqrt{2}}  & \frac{q\cos{(2k+\nu)}}{2}& \frac{\sqrt{6}q\sin{(2k+\nu)}}{4}  \\
0 & -\frac{q\sin(2k+\nu)}{2}&-\frac{q\cos{(2k+\nu)}  }{\sqrt{2}}  &\frac{3q\sin{(2k+\nu)}  }{4}  &0&  0& -\frac{q\cos{(2k+\nu)}  }{\sqrt{2}} &-\frac{3q\sin{(2k+\nu)}  }{4\sqrt{12}} . \\
\end{pmatrix}}.\]\\

Next we will show that the superoperator $\mathcal{L}_{k,k}$ satisfies the eigenvalue condition, it's sufficient to show $\mathcal{L}_{k,k,1}$ satisfies the eigenvalue condition. We will use the standard basis $|\xi_{i}><\xi_j|,  i,j=1, 2 ,\dots, n,$ with lexicographical ordering \[|\xi_{1}><\xi_1|,  |\xi_{1}><\xi_2, \dots, |\xi_{1}><\xi_n|, |\xi_{2}><\xi_1|, \dots, |\xi_{2}><\xi_n|, \dots, |\xi_{n}><\xi_1|, \dots, |\xi_{n}><\xi_n|.\] \\

Note that  $\mathcal{L}_{k,k,1} (\rho)=B_1\rho {B_1}^*+B_2\rho {B_2}^*= \left(B_1\otimes \overline{B_1}+ B_2 \otimes \overline{B_2}\right) (\rho)$. Then 
\[ \mathcal{L}_{k,k,1} (\rho)= r^2 \begin{pmatrix}
 1 &0   &1&0&0 &0 &1 & 0&1\\
  r & 1   &-r&0&0 &0 &r & 1&-r\\
 0 &0   &0&0&0 &0 &0 & 0&0\\
r &0   &r&1&0 &1 &-r & 0&-r\\
  r^2 &r   &-r^2&r&1 &-r &-r ^2&-r & r^2\\
 0 &0   &0&0&0 &0 &0 & 0&0\\
0 &0   &0&0&0 &0 &0 & 0&0\\
0 &0   &0&0&0 &0 &0 & 0&0\\
 r^2 & -r   &-r^2&-r&1 &r &-r ^2&r & r^2\\
\end{pmatrix}.\]

 Let $\lambda$ be an eigenvalue of $\mathcal{L}_{k,k,1}$ and $\beta=\frac{\lambda}{r^2}$. Then 
 \[ \frac{\mathcal{L}_{k,k,1}}{r^2}-\beta= \begin{pmatrix}
 1-\beta &0   &1&0&0 &0 &1 & 0&1\\
  r & 1-\beta   &-r&0&0 &0 &r & 1&-r\\
 0 &0   &-\beta&0&0 &0 &0 & 0&0\\
r &0   &r&1-\beta&0 &1 &-r & 0&-r\\
  r^2 &r   &-r^2&r&1-\beta &-r &-r ^2&-r & r^2\\
 0 &0   &0&0&0 &-\beta &0 & 0&0\\
0 &0   &0&0&0 &0  & -\beta&0&0\\
0 &0   &0&0&0 &0 &0 & -\beta&0\\
 r^2 & -r   &-r^2&-r&1 &r &-r ^2&r & r^2-\beta\\
\end{pmatrix}.\]

By using the minors matrices on rows with only nonzero entry $-\beta$, we have 
\[\frac{\mathcal{L}_{k,k,1}}{r^2}-\beta=
\beta^4 \left |\begin{matrix} 
 1-\beta &0   &0 &0 &1 \\
  r & 1-\beta  & 0&0 &-r\\
r &0  &1-\beta&0 &-r\\
  r^2 &r  &r&1-\beta & r^2\\
 r^2 & -r  &-r&1 & r^2-\beta\\
\end{matrix} \right |=\frac{1}{2}(\beta-1)(\beta-2)(-2\beta^2+5\beta-4)(\beta+1).\]
Hence $\beta=2$ is the largest simple root. Equivalently, $\lambda=1$ is the largest eigenvalue of $\mathcal{L}_{k,k,1}$ with multiplicity one and the absolute value of all other eigenvalue are strictly less than 1. That is, $\mathcal{L}_{k,k,1}$ satisfies the eigenvalue condition.
By  Theorem (\ref{eigencond1}), the superoperator $\mathcal{L}_{k,k}$ satisfies the eigenvalue condition for all  $0<p\le 1$ and by Theorem \ref{CONV}, we have 
\[
\hat{P}(\frac{\nu}{\sqrt{t}},t)\rightarrow\frac{1}{2\pi}\int_0^{2\pi}e^{-\frac{1}{2}\sigma^2 \nu^2}dk,
\]
where $\sigma^2=\frac{1}{I(1,0)}\left [\frac{13}{3}I(1,0)-\frac{4}{9}I_\lambda(1,0)-4\sqrt{2}J(1,0)+\frac{4\sqrt{2}}{9}J_\lambda(1,0)-\frac{4\sqrt{2}i}{3}J_\nu(1,0)-\frac{1}{3}\right ].$

\end{example}


\begin{Ack} Wei-Shih Yang would like to  thank the Institute
of Mathematics, Academia Sinica, for its support during his visit, while part of this work was done.\end{Ack}

\end{document}